\documentclass[fleqn]{mat01}
\usepackage{times,mathtimy,amssymb,epsfig}
\begin{document}

\setcounter{page}{153}
\firstpage{153}

\newtheorem{theore}{Theorem}
\renewcommand\thetheore{\arabic{section}.\arabic{theore}}
\newtheorem{theor}[theore]{\bf Theorem}
\newtheorem{propo}[theore]{\rm PROPOSITION}
\newtheorem{lem}[theore]{Lemma}
\newtheorem{definit}[theore]{\rm DEFINITION}
\newtheorem{coro}[theore]{\rm COROLLARY}
\newtheorem{rem}[theore]{Remark}
\newtheorem{exampl}[theore]{Example}

\newcommand{\Z}{\mathbb{Z}}
\newcommand{\R}{\mathbb{R}}
\newcommand{\C}{\mathbb{C}}
\newcommand{\N}{\mathbb{N}}
\renewcommand{\L}{\mathcal{L}}
\renewcommand{\S}{\mathcal{S}}
\renewcommand{\H}{{\mathbb{H}}}
\newcommand{\F}{\mathfrak{F}}
\newcommand{\Zp}{\mathbb{Z}/p\mathbb{Z}}
\newcommand{\Zq}{\mathbb{Z}/q\mathbb{Z}}
\newcommand{\Zz}{\mathbb{Z}/2\mathbb{Z}}
\newcommand{\Zn}{\mathbb{Z}/n\mathbb{Z}}
\newcommand{\sm}{\setminus}
\newcommand{\tr}{\triangle}
\newcommand{\co}{\colon\thinspace}
\newcommand{\del}{\partial}

\font\zz=msam10 at 10pt
\def\Box{\mbox{\zz{\char'244}}}

\font\bi=tibi at 10.4pt

\title{Limits of functions and elliptic operators}

\markboth{Siddhartha Gadgil}{Limits of functions and elliptic operators}

\author{SIDDHARTHA GADGIL}

\address{Stat Math Unit, Indian Statistical Institute,
Bangalore~560~059, India\\
\noindent E-mail: gadgil@isibang.ac.in}

\volume{114}

\mon{May}

\parts{2}

\Date{MS received 10 February 2004}

\begin{abstract}
We show that a subspace $S$ of the space of real analytical functions on
a manifold that satisfies certain regularity properties is contained in
the set of solutions of a linear elliptic differential equation. The
regularity properties are that $S$ is closed in $L^2(M)$ and that if a
sequence of functions $f_n$ in $S$ converges in $L^2(M)$, then so do the
partial derivatives of the functions $f_n$.
\end{abstract}


\maketitle

\vspace{.3pc}
\noindent \looseness -1 The limit $f$ of a\vspace{-.05pc} sequence $f_n$ of complex analytical functions (under
uniform convergence on compact sets) is complex analytical. Furthermore
all partial derivatives of $f_n$ converge to the corresponding partial
derivatives of $f$. This is in contrast to the case of real analytical
functions. In fact, by the Weierstrass approximation theorem, every
continuous real function on a compact domain is the uniform limit of
real analytical functions on this domain.

The reason for the contrast\vspace{-.05pc} between the complex and the real analytical
cases is of course that complex analytical functions satisfy an elliptic
differential equation, namely the Cauchy--Riemann equation (or
alternatively because they satisfy the Laplace equation), while no such
equation is satisfied in the analytical case.

Here we show that this phenomenon\vspace{-.05pc} is universal, namely, whenever we have
a class $S$ of (real analytical) functions on a closed manifold $M$ that
have regularity properties similar to those of holomorphic functions,
all functions in $f\in S$ satisfy an elliptic differential equation
$Pf=0$.

Our motivation is that in many\vspace{-.05pc} geometric situations rigidity phenomena
are associated with elliptic operators which are often \emph{hidden},
i.e., not \emph{a priori} related to the geometry. Two striking
instances of this are the Seiberg--Witten equations for smooth
four-dimensional manifolds and $J$-holomorphic curves in Symplectic
topology. Hence it is of interest to show that there are situations
where there must be elliptic operators, even though they are not \emph{a
priori} present.

First we recall the definition of the Sobolev spaces $W^{2,k}$ where
$k\geq 0$ is an integer. We will not need the general case when $k\in
R$.

\begin{definit}$\left.\right.$\vspace{.5pc}

\noindent {\rm Let $k\geq 0$ be an integer. Suppose $f$ and $g$ are
smooth, real valued functions on $\R^n$ with compact support, we define
the Sobolev inner product $\langle f,g\rangle_{2,k}$ by\vspace{-.2pc}
\begin{equation*}
\langle f,g\rangle_{2,k} = \Sigma_{j=0}^k \Sigma_{|I|=j}\int_{\R^n}
\del^If(x) \del^Ig(x)\,{\rm d}x,
\end{equation*}

$\left.\right.$\vspace{-.7pc}

\noindent where $I$ is a multi-index and $\del^I$ denotes the partial derivative
with respect to $I$.}\vspace{-.75pc}
\end{definit}

\pagebreak

\begin{definit}$\left.\right.$\vspace{.5pc}

\noindent {\rm Suppose $M$ is a manifold, let $\{U_i\}$ be a locally
finite cover of $M$ by subsets homeomorphic to $\R^n$ and let
$\{\pi_i\}$ be a partition of unity subordinate to this cover with
${\rm supp}(\pi_i)\subset U_i$ compact. For smooth compactly supported
functions $f$ and $g$ on $M$, define $\langle f,g\rangle_{2,k}$ by
\begin{equation*}
\langle f,g\rangle_{2,k}=\Sigma_i \langle \pi_i f,\pi_i g\rangle_{2,k},
\end{equation*}
where $\langle \pi_i f,\pi_i g\rangle_{2,k}$ denotes the Sobolev inner
product on $U_i=\R^n$.}
\end{definit}

The above definition depends on the choice of the cover $U$, but
different covers give equivalent inner products.

\begin{definit}$\left.\right.$\vspace{.5pc}

\noindent {\rm The Sobolev space $W^{2,k}(M)$ is the Hilbert space
completion of the space $C^\infty_c(M)$ of smooth functions on $M$ with
compact support with respect to the Sobolev inner product $\langle \ ,\
\rangle_{2,k}$.}
\end{definit}

When $k=0$ we get the Hilbert space $L^2(M)$ with its usual inner
product. The definitions above coincide with the definitions using
Fourier transforms.

We can now state our main result.

\begin{theor}[\!]
Let $S$ be a subspace of real analytical functions on a compact real
analytical manifold $M$ that is closed under the $L^2$-norm on
$M$. Assume further that if $f_n\in S$ is a sequence of functions such
that $f_n\to f$ in $L^2(M)${\rm ,} then $f_n\to f$ in all Sobolev spaces
$W^{2,k}(M)${\rm ,} $k\in \N$. Then there is an analytical elliptic
differential operator $P$ on $M$ such that $\forall f\in S, Pf=0$.
\end{theor}

\begin{rem}{\rm
The analogous result for sections of a bundle on $M$ holds and can be
proved in exactly the same way.}
\end{rem}

A differential operator $P$ that satisfies elliptic regularity on
every open set $U$ (i.e., if $u$ is a distribution on $U$ with $Pu=f$,
$f$ smooth, then $u$ is smooth) is called \emph{hypoelliptic}. Such
operators have been characterised among operators with constant
coefficients by H\"ormander~\cite{Ho}. What we consider here is a
different situation where our class of functions may not be given by a
differential equation. What we can conclude is also weaker -- we only
know that $S$ is contained in the set of solutions to an elliptic
differential equation.

We now outline the proof. By using the hypothesis, we show that on the
space $S$, the $L^2$ norm is equivalent to the $W^{2,2}$ norm. From
this we deduce that the space $S$ is finite dimensional. Next, for
each $x\in M$, the partial derivatives at $x$ give linear functionals
on $S$. By using the finite-dimensionality of $S$, we show that at $x$
we can find an elliptic differential equation satisfied by $S$. The
same method yields elliptic differential equations on certain
semi-analytical sub-varieties. Finally, we use the local Noetherian
property of real analytic varieties to deduce that we can globally
construct an elliptic differential operator $P$ with $Pf=0\ \forall
f\in S$.

Only the final step in the above outline uses analyticity. We shall
show, however, that the hypothesis of analyticity is essential for our
result.

\section{Finite dimensionality of {\bi S}}

In this section we show that $S$ is finite dimensional. First we make
an elementary observation about subspaces of Hilbert spaces.

Let $H_1$ and $H_2$ be Hilbert spaces with norms $\|\ \|_1$ and $\|\
\|_2$ respectively. Assume that as sets $H_2\subset H_1$ with
$\|x\|_1\leq \|x\|_2$, $\forall x\in H_2$. The following result will be
applied to the case when $H_1=L^2(M)$ and $H_2=W^{2,2}(M)$.

\setcounter{theore}{0}
\begin{propo}$\left.\right.$\vspace{.5pc}

\noindent Let $S$ be a subspace of $H_2$ that is closed in $H_1$ and
$H_2$ so that the subspace topologies induced by $H_1$ and $H_2$
coincide. Then there is a constant $C>0$ such that $\|x\|_2\leq
C\|x\|_1, \forall x\in S$.
\end{propo}

\begin{proof}
By hypothesis, the identity map from $S$ with its topology as a subspace
of $H_1$ to $S$ with its topology as a subspace of $H_2$ is continuous.
Hence it must be bounded from which the conclusion follows.\hfill $\Box$
\end{proof}

Note that by the hypothesis in our main theorem, the above result
applies to $S$ with $H_1=L^2(M)$ and $H_2=W^{2,2}(M)$. We next show that
a space $S$ satisfying the hypothesis of the main theorem is finite
dimensional.

\begin{lem}
Let $M$ be a closed manifold and let $S$ be a subspace of
$W^{2,2}(M)\subset L^2(M)$ such that there exists $C>0$ such that for
$f\in S, \|f\|_{2,2}<C\|f\|_2$. Then $S$ is finite dimensional.
\end{lem}

\begin{proof}
Suppose $S$ is infinite dimensional, then let $\{x_n\}$ denote an
$L^2$-orthonormal sequence in $S$. By hypothesis, for all $j\in \N$,
$x_n\in W^{2,2}$ and $\|x_n\|_{2,2}<C$. By the Rellich lemma it follows
that the sequence $x_n$ has a convergent sequence in $L^2$. But, as the
vectors $x_n$ are $L^2$-orthonormal, this is impossible. Thus, $S$ must
be finite dimensional.\hfill $\Box$
\end{proof}

\section{Pointwise differential equations}

We can now construct elliptic differential equations satisfied by the
functions in $S$ at a single point $x\in M$. Choose a system of local
coordinates. Observe that partial derivatives at $x$ give linear
functionals on $S$, i.e., elements of the dual $S'$ of $S$. These
generate a subspace $V_x$ of $S'$. As $S'$ is finite dimensional, $V_x$
is generated by finitely many partial derivatives, and hence those of
order at most $k$ for some $k$. We denote these differential operators
by $P_1$,\dots, $P_m$.

Now let $E$ be an elliptic operator of order greater than $k$, for
instance a power of the Laplacian. Then $f\mapsto Ef(x)$ is an element
of $S'$, hence is spanned by $P_i$. Thus, at $x$ each $f\in S$
satisfies a relation $(E-\Sigma c_iP_i)f(x)=0$. As this has the same
leading term as $E$, this is an elliptic differential equation.

Note that the relations $P_1$,\dots, $P_m$ are independent as elements
of $S'$ on an open set (as independence is an open condition). Let $r(x)
= {\rm dim} V_x$. Let $f_1$,\dots,$f_N$ be a basis for $S$. In the
special case where $r(x)$ is a constant function (for instance $r(x) =
{\rm dim}(S)$, the maximum possible value, at all points), we shall see
that we have a global elliptic operator even in the absence of
analyticity.

\setcounter{theore}{0}
\begin{propo}$\left.\right.$\vspace{.5pc}

\noindent Suppose $r(x)=m$ is a constant. Then there is an elliptic
differential operator $E$ such that $Ef=0${\rm ,} for all $f\in S$.
\end{propo}

\begin{proof}
We first show that there is a uniform degree $k$ so that the operators
of degree at most $k$ span $V_x$ for all $x\in M$. For $x\in M$, let
$P_1$,\dots, $P_m$ be operators independent at $x$ and let $U_x$ be the
set where these operators are independent. This is an open set as linear
independence is an open condition (for instance, by considering
determinants). By hypothesis, $M$ is the union of the sets $U_x$. By
compactness we can find finitely many such sets $U_j$ whose union is
$M$. Let $k$ be the maximal degree of the differential operators
associated to these sets.

Now, let $E$ be an elliptic operator of order greater than $k$. For each
$U_j$, we have differential operators $P_1$,\dots,$P_m$ which are
independent at each $x\in U_j$ and hence span $V_x$. Hence we have a
relation $Ef(x)=\Sigma c_i(x) P_i f(x)$ for $f\in S$.

We next show that each $c_i(x)$ is smooth as a function of $x\in U_j$.
Let $x_0\in U_j$ be an arbitrary point. We shall show that $c_i(x)$ is
smooth at $x_0$.

As the operators $P_j$, $1\leq j\leq m$, are independent at $x_0$ as
functionals on $S$, there exist $g_i\in S$, $1\leq i\leq m$, such that
$P_jg_i(x_0)=\delta_{ij}$. It follows that for $x\in U_j$, $P_j
g_i(x)=\delta_{ij}+a_{ij}(x)$ with $a_{ij}(x)$ smooth functions and
$a_{ij}(x_0)=0$.

Let $A(x)$ denote the matrix with entries $a_{ij}(x)$ and let $V(x)$
(respectively $C(x)$) denote the (column) vector with entries $Eg_i(x)$
(respectively $c_i(x)$). Note that $V(x)$ is smooth as a function of
$x$. As $Eg_i(x)=\Sigma_j P_j g_i(x) c_j(x)$, i.e., $V(x)=(I+A(x))C(x)$,
we have $C(x)=(I+A(x))^{-1}V(x)$.

Now, $A(x_0)=0$ and it is well-known that $M\mapsto (I+M)^{-1}$ is
smooth as a function of $M$ at $M=0$ (by using the implicit function
theorem or Taylor expansions). Hence $C(x)$ is smooth at $x_0$, i.e.,
each $c_i$ is smooth at $x_0$, as required. Let $E_j'=E-\Sigma_i c_i
P_i$. This is an elliptic operator annihilating $S$.

To construct an elliptic operator globally, we take a partition of unity
$\{\pi_i\}$ subordinate to the cover $\{U_j\}$ and let
$E_0=\Sigma_i\pi_iE_i$. Then each $f\in S$ is in the kernel of $E_0$ and
$E_0$ is elliptic as, by construction, the leading term of $E_0$ is the
same as that of $E$.\hfill $\Box$
\end{proof}

Without assuming analyticity, however, our main result fails in general.
To see this, we let $M=S^1=\R/\Z$ and construct a function $f$ on $S^1$
so that, for $n>1$, $f^{(n)}(1/n)\neq 0$ but $f^{(k)}(1/n)=0$ for $k<n$.
Let $S$ be the (one-dimensional) span of $f$.

An elliptic differential operator $E$ on a one-dimensional manifold is a
differential operator $P(D)$ whose leading coefficient is non-zero at
all points. The function $f$ does not satisfy $Ef=0$ for any such
operator as if $d$ is the order of $E$, by construction $Ef(1/d)\neq 0$.

\section{Globalisation in the analytical case}

In the analytical case, we shall construct sets similar to $U_j$ above.
These are now open in the real analytic topology, i.e., one whose
sub-basis is generated by sets of the form $f(x)\neq 0$ where $f$ is an
analytical function.

We shall need two basic facts regarding the real analytical topology
(see, for instance \cite{Na}). Firstly, any closed set is defined by a
single equation, as given a closed set $F = Z(g_1, \dots, g_p) =
\{x:g_i(x) = 0, 1\leq i\leq p\}$, we have $F = Z(g_1^2 + \dots +
g_p^2)$. Secondly, as the ring of power series is Noetherian, the real
analytical topology is locally Noetherian. As $M$ is compact, the real
analytical topology on $M$ is Noetherian.

As $M$ is analytical, by a theorem of Morrey~\cite{Mo} and
Grauert~\cite{Gr} there is an analytical Riemannian metric on $M$. Hence
the Laplacian (with respect to an analytical metric) is an analytical
elliptic operator on $M$ and so are it powers. It follows that there are
analytical elliptic differential operators on $M$ of arbitrarily high
orders.

In the rest of this section, we make the convention that all
differential operators we consider are analytical ones globally defined
on $M$. In particular we shall use the notation of the previous
sections, but with $V_x$ now being the subspace of $S'$ generated by
global analytical operators acting on $S$ at $x$ and $r(x)$ its
dimension.

We shall inductively construct sequences of sets $F_i$ and $V_i$, with
$F_i$ a decreasing sequence of closed sets and $V_i$ open, and finitely
many elliptic differential operators that span $V_x$ for $x\in F_i\cap
V_i$.

Let $F_1=M$ and note that this is a closed subset of $M$. On $F_1=M$,
let $m_1=r(z)$ be the maximum value of $r(x)$ (which is attained as
$r(x)\in\Z$, $0\leq r(x)\leq {\rm dim}(S)$ for all $x\in M$) and let
$P^1_1$,\dots, $P^1_m$, $m=m_1$, be (analytic) differential operators
with $f\mapsto P^1_if(z)$ independent in $S'$. Then the set $V_1$ where
the $P^1_j$'s are independent is an open set in the analytical topology.
It follows that for all $x\in V_1$, the functionals $P^1_j(x)$ span
$V_x$. Let $q_1$ be the maximum order of differential operators $P^1_j$.

Next, let $F_2$ be the complement of $V_1$. Let $r(z)=m_2\leq m_1$,
$z\in F_2$, be the maximum value of $r(x)$ on $F_2$ and let
$P^2_1,\dots,P^2_m, m=m_2$, be different operators with $f\mapsto
P^2_jf(z)$ independent functionals in $S'$. Let $V_2$ be the open set
(in the real analytical topology) where $P^2_j$'s are independent as
functionals in $S'$. Then these span $V_x$ for all $x\in F_2\cap V_2$.
Let $q_2$ be the maximum order of these differential operators.

Inductively, given $F_k$ and $V_k$, we define $F_{k+1} = F_k\!\sm\! V_k$
and let $m_{k+1}$ be the maximum rank of $V_x$ on $F_{k+1}$. As above,
we construct differential operators $P^{k+1}_j$ and let $V_{k+1}$ be the
set on which these are independent. These span $V_x$ for all $x\in
F_{k+1}\cap V_{k+1}$. Let $q_k$ be the maximal order of the differential
operators constructed above.

Now, by the local Noetherian property, the above process must stabilise.
It follows that for some $n$, $F_n\subset V_n$. Let $q$ be the maximum
of the numbers $q_j$, $1\leq j\leq n$ and note that on each set $F_j\cap
V_j$, we have differential operators of degree at most $q$ that span at
each point the subspace $V_x$. Let $g_n$ be analytical functions such
that $F_n=Z(g_n)$.

Let $E$ be an analytic elliptic operator of order greater than $q$. We
construct inductively analytic elliptic operators $E_n$, \dots, $E_1$
with $E_1$ being the required operator. First, note that on $F_n\cap
V_n$, we can find an operator $G_n$, with analytical coefficients, of
order at most $q$ so that $Ef(x)=G_nf(x)$ for all $x\in F_n\cap V_n$,
$f\in S$. Let $E_n=E-G_n$.

Next, observe that on $F_{n-1}\cap V_{n-1}$, the function $g_n$ does not
vanish, and hence the operator $E_n/g_n$ is well-defined. Hence there is
an operator $G_{n-1}$, with analytical coefficients, of order at most
$k$ so that $(E_n/g_n)f(x)=G_{n-1}f(x)$ for all $x\in F_{n-1}\cap
V_{n-1}$, $f\in S$. Let $E_{n-1}=E_n-g_nG_{n-1}$. This annihilates $S$
on $F_{n-1}\cap V_{n-1}$ by construction and also on $F_n$ as it
coincides with $E_n$ on $F_n$. Thus $E_{n-1}$ annihilates $S$ on
$F_{n-1}$.

Similarly, given an elliptic operator $E_k$ that annihilates $S$ on
$F_k$, we can construct an elliptic operator $E_{k-1}$ that annihilates
$S$ on $F_{k-1}$. Proceeding inductively, we obtain an operator that
annihilates $S$ on $F_1=M$.\hfill $\Box$

\section*{Acknowledgements}

The author would like to thank V~Pati and Harish Seshadri for numerous
helpful comments and conversations.

\end{document}